\documentclass[11pt]{amsart}
\usepackage{amsmath,amssymb,amsthm,color}
\usepackage[T1]{fontenc}
\newtheorem{theorem}{Theorem}[section]
\newtheorem{lemma}[theorem]{Lemma}
\newtheorem{proposition}[theorem]{Proposition}
\theoremstyle{remark}
\newtheorem{remark}[theorem]{Remark}
\newtheorem{example}[theorem]{Example}
\newcommand{\C}{\mathbb C}
\DeclareMathOperator{\Fl}{Fl}
\DeclareMathOperator{\Hess}{Hess}
\DeclareMathOperator{\Aut}{Aut}
\DeclareMathOperator{\Lie}{Lie}
\DeclareMathOperator{\wt}{wt}
\DeclareMathOperator{\spanop}{span}

\newcommand{\e}{\mathbf e}
\newcommand{\U}{\mathbb U}

\title[Automorphisms of Hessenberg varieties]
{Global vector fields and the automorphism group of regular semisimple Hessenberg varieties}
\author{Mikiya Masuda}
\author{Takashi Sato}
\author{Haozhi Zeng}

\begin{document}

\begin{abstract}
For a connected regular semisimple Hessenberg variety other than the full flag variety, we prove that the identity component of its automorphism group is an algebraic torus.
The argument uses equivariant local coordinates and weight-homogeneous coefficient functions of vector fields.
\end{abstract}
\maketitle

\section{Introduction}
Hessenberg varieties form a class of subvarieties of flag varieties defined by linear incidence conditions.
Their name comes from Hessenberg matrices, which play a basic role in numerical linear algebra,
but their geometry brings together questions from algebraic geometry, representation theory, and algebraic combinatorics.
They were introduced and systematically studied by De Mari, Procesi, and Shayman \cite{DMPS}.
In type A, a Hessenberg variety is determined
by an endomorphism $S$ of $\C^n$ and a weakly increasing function
$h:[n] =\{1,\dots,n\}\to[n]$ satisfying $h(i)\ge i$ for any $i\in [n]$;
it consists of complete flags $F_\bullet$ such that
\[
	SF_i\subseteq F_{h(i)}\qquad(1\le i\le n).
\]
This single construction contains several familiar spaces.
For example, $h(i)=n$ for any $i\in [n]$ gives the full flag variety $\Fl(\C^n)$,
whereas $h(i)=i$ for any $i\in [n]$ gives the variety of $S$-stable flags.

The geometry depends strongly on the Jordan type of $S$.
The regular nilpotent case is closely related to Springer theory and includes the Peterson variety.
The present paper concerns the opposite generic case, in which $S$ is regular semisimple.
In this case $S$ may be taken to be diagonal with pairwise distinct eigenvalues, and $\Hess(S,h)$ is a smooth projective variety \cite{DMPS}.
Its topology admits a particularly effective combinatorial description.
Affine pavings for Hessenberg varieties were developed by Tymoczko \cite{TymoczkoLinear},
and the diagonal torus action makes a regular semisimple Hessenberg variety a GKM space.
Thus its equivariant cohomology can be read from a labelled graph in
the sense of Goresky--Kottwitz--MacPherson \cite{GKM}.
Tymoczko's dot action then equips its cohomology with an action of the symmetric group \cite{TymoczkoDot}.
This action led to the connection with chromatic quasisymmetric functions discovered by Shareshian and Wachs \cite{SW};
the corresponding conjecture was proved by Brosnan and Chow \cite{BC}.
These results illustrate how the torus action and the combinatorics of the Hessenberg function control much of the topology of the variety.

It is therefore natural to ask how much of this combinatorial symmetry is
realized by algebraic automorphisms.  Let
\[
	X = \Hess(S,h),\qquad G = \Aut^0(X),
\]
where $\Aut^0(X)$ denotes the identity component of the automorphism group $\Aut(X)$ of $X$,
and assume that $X$ is connected and is not the full flag variety.
The quotient of the diagonal torus
\[
 	\mathbb T = (\C^*)^n/\C^*
\]
acts effectively on $X$, so $\mathbb T \subset G$.
In \cite{JangEtAl}, automorphisms of the GKM graph were used to study $\Aut(X)$.
In particular, it was shown that the reductive part of $G$ is $\mathbb T$ and that
\[
	G \cong \mathbb T \ltimes G_u,
\]
where $G_u$ is the unipotent radical.
The remaining issue is whether $G_u$ can be nontrivial.
This is not a formal consequence of the description of the reductive part:
an algebraic group with reductive part $\mathbb T$ may still have a nontrivial unipotent radical.
Thus the determination of $G_u$ is precisely the part of the identity component left open by \cite{JangEtAl}.
The purpose of this paper is to settle this remaining problem for every connected regular
semisimple Hessenberg variety of type A other than the full flag variety, and
thereby determine the identity component completely:
\[
	\Aut^0(\Hess(S,h)) = \mathbb T \cong (\C^*)^{n-1}.
\]

The codimension-one case has also been studied by Brosnan, Escobar, Hong,
Lee, Lee, Mellit, and Sommers \cite{BEHLLMS} in their broader investigation of
automorphisms and deformations of regular semisimple Hessenberg varieties.
They treat arbitrary simple Lie types and, in addition to determining the
identity component of the automorphism group, compute the first-order
deformation space and prove the vanishing of the higher tangent cohomology.
In type A, they determine the full automorphism group and describe the corresponding moduli stack.
Their identity-component result agrees with the theorem above in the codimension-one case, whereas the present argument applies to every connected regular semisimple Hessenberg variety of type A other than the full flag variety.

The proof is local and algebraic.
Working on the opposite big cell around the identity fixed point, we represent a flag by a lower unitriangular matrix $g = (g_{kl})$.
The matrix $A=g^{-1}Sg$ represents the action of $S$ with respect to the basis of $\C^n$ consisting of the column vectors of $g$, and the Hessenberg condition simply translates to $A_{kl}=0$ for $k>h(l)$.
The equation $Sg=gA$ then yields a system of equations that can be solved recursively by induction on the height $k-l$.
This shows directly that the entries $g_{kl}$ for $l < k \leq h(l)$ serve as coordinates on the chart.
Furthermore, for each such index $(i, j)$, we consider the subtorus $\mathbb{T}_{ij} \subset \mathbb{T}$
whose entries from $j$ to $i$ are equal.
The fixed-point locus of $\mathbb{T}_{ij}$ contains a component naturally isomorphic to the smaller flag variety $\Fl(\spanop(\e_j,\dots,\e_i))$, and restricting to this component simply sets the coordinates $g_{kl}$ outside this block to zero.

These coordinates allow us to detect invariant vector fields on $X$.
A nonzero $\mathbb T$-invariant global vector field has a nonzero restriction to at least one of
the smaller flag varieties above.
On the other hand, an additive one-parameter subgroup centralizing $\mathbb T$ must act trivially on every such flag variety.
This excludes additive subgroups from $C_G(\mathbb T)\cap G_u$, where $C_G(\mathbb T)$ denotes the centralizer of $\mathbb T$ in $G$.
A weight argument at the torus fixed points then gives $C_G(\mathbb T)=G$ unless $X$ is the full flag variety.
Since every nontrivial connected unipotent group over $\C$ contains a copy of
the additive group $\mathbb G_a \cong \C$, it follows that $G_u$ is trivial.

The paper is organized as follows.
Section~2 sets up the notation for regular semisimple Hessenberg varieties and their vector fields.
Section~3 constructs the affine coordinates and the smaller flag varieties.
Section~4 uses these coordinates to rule out centralizing additive subgroups.
Section~5 computes the centralizer of the diagonal torus,
and Section~6 proves the main theorem.

\section{Setup and vector fields}
Let $h \colon [n] \to [n]$ be weakly increasing with $h(i)\ge i$ for all $i \in [n]$,
and let $S=\operatorname{diag}(\lambda_1,\ldots,\lambda_n)$ have pairwise distinct eigenvalues. Put
\[
	X = \Hess(S,h)
	= \{F_\bullet\in\Fl(\C^n) \mid SF_i \subseteq F_{h(i)} \text{ for all } i\}.
\]
Assume $h(i)\ge i+1$ for $i<n$, so that $X$ is connected.
Since $X$ is smooth, connectedness implies irreducibility.
The effective diagonal torus $\mathbb T=(\C^*)^n/\C^*$ acts on $X$.
We write $G=\Aut^0(X)$ and $G_u$ for its unipotent radical.
We denote by $t_i$ the standard characters of $(\mathbb{C}^*)^n$, so that the differences $t_i - t_j$ are in the character group of $\mathbb T$.

The fixed points of $\mathbb T$ are permutation flags indexed by permutations $w\in S_n$.
Set
\[
	H = \{(i,j) \mid 1\le j<i\le h(j)\}.
\]
At the identity fixed point $e$, the tangent representation is
\[
	T_{e}X = \bigoplus_{(i,j)\in H} \C_{t_i-t_j},
\]
where $\C_{t_i-t_j}$ denotes the complex $\mathbb T$-module with weight $t_i-t_j$. 
The infinitesimal action identifies the Lie algebra $\Lie(G)$ of $G$ with $H^0(X,TX)$, the space of global algebraic vector fields of $X$, equivariantly for conjugation by $\mathbb T$.
We will also use the fact that a nonzero global vector field on an irreducible variety has a nonzero germ at every point:
if its germ vanishes at one point, it vanishes on a nonempty open subset and hence everywhere.

\section{Smaller flag varieties and equivariant coordinates}
For $(i,j)\in H$, let $V_{ij}=\spanop(\e_j,\ldots,\e_i)$,
where $\e_1,\dots,\e_n$ denote the standard basis of $\C^n$,
and let $\mathbb T_{ij}\subset \mathbb T$ be the subtorus on which the entries indexed by $j,\ldots,i$ are equal.
Let $X_{ij}$ denote the connected component of $X^{\mathbb T_{ij}}$ containing the identity fixed point $e$.

\begin{lemma}\label{lem:smaller}
For every $(i,j)\in H$, the variety $X_{ij}$ is naturally isomorphic to $\Fl(V_{ij})$. Moreover,
\[
	T_eX_{ij}=\bigoplus_{j\le l<k\le i}\C_{t_k-t_l}.
\]
\end{lemma}
\begin{proof}
The connected component of $\Fl(\mathbb{C}^n)^{\mathbb{T}_{ij}}$ containing the standard flag $e$ consists precisely of the flags
\[
	F_k =
	\begin{cases}
		\spanop(\e_1,\ldots,\e_k),	&k<j,\\
		\spanop(\e_1,\ldots,\e_{j-1})\oplus E_{k-j+1},	&j\le k \le i,\\
		\spanop(\e_1,\ldots,\e_k),	&k>i,
	\end{cases}
\]
where $E_\bullet$ is an arbitrary complete flag in $V_{ij}$.
This component is isomorphic to $\Fl(V_{ij})$.
For $j\le k\le i$, monotonicity gives $h(k)\ge h(j)\ge i$, and hence $SF_k\subseteq F_i\subseteq F_{h(k)}$; outside this interval, $F_k$ is $S$-stable.
Thus all these flags belong to $X$, and hence this component coincides with $X_{ij}$.
Consequently, the tangent weights at $e$ are precisely $t_k-t_l$ for $j\le l<k\le i$.
\end{proof}

We now construct an algebraic affine chart. Let $U\subset\Fl(\C^n)$ be the opposite big cell containing $e$.
A flag in $U$ is represented uniquely by a lower unitriangular matrix $g=(g_{kl})$, with $g_{kk}=1$ and $g_{kl}=0$ for $k<l$. Write $A=g^{-1}Sg$.
The Hessenberg conditions are precisely
\begin{equation}\label{eq:hess}
	A_{kl}=0\qquad(k>h(l)).
\end{equation}
Since $g$ is lower unitriangular and $S$ is diagonal, $A$ is lower triangular with diagonal entries $A_{kk}=\lambda_k$.
Comparing the $(k,l)$ entries of $Sg=gA$ for $k>l$ gives
\begin{equation}\label{eq:recursion}
	(\lambda_k-\lambda_l)g_{kl}=A_{kl}+\sum_{l<m<k}g_{km}A_{ml}.
\end{equation}

We call $k-l$ the height of the $(k,l)$-entry of a matrix.
Note that, when $(k,l) \in H$, we have $(k,m), (m,l) \in H$ for $l<m<k$.
When $(k,l) \notin H$ with $k > l$, recall that $A_{kl} = 0$ by \eqref{eq:hess}.
Hence, using \eqref{eq:recursion}, we can express the non-zero entries of $A$ as a polynomial in $\{g_{kl}\mid (k,l) \in H\}$ and entries of $A$ of strictly smaller height.
Inductively, every entry of $A$ can be expressed as a polynomial in $\{g_{kl}\mid (k,l) \in H\}$.

On the other hand, for $(k,l) \notin H$ with $k > l$,
we can express $g_{kl}$ as a polynomial in entries of $g$ and $A$ of strictly smaller height, using \eqref{eq:recursion} again.
Inductively, combined with the expressions for $A$, every entry of $g$ can be expressed as a polynomial in $\{g_{kl}\mid (k,l) \in H\}$.

\begin{example}
	Let $h = (3,3,4,4)$, where we write $h$ as the sequence of its values $(h(1), h(2), \ldots, h(n))$.
	In this case, the argument above expresses $g_{41}$ and $g_{42}$ as polynomials in $\{g_{kl} \mid (k,l) \in H\}$ as follows.
	First, we express the relevant entries of $A$ as polynomials in $\{g_{kl} \mid (k,l) \in H\}$, using \eqref{eq:recursion}.
	Concretely, $A_{21} = (\lambda_2 - \lambda_1)g_{21}$,
	$A_{32} = (\lambda_3 - \lambda_2)g_{32}$, $A_{43} = (\lambda_4 - \lambda_3)g_{43}$, and 
	\[
		A_{31} = (\lambda_3 - \lambda_1)g_{31} - g_{32}A_{21}
		= (\lambda_3 - \lambda_1)g_{31} - (\lambda_2 - \lambda_1)g_{32}g_{21}.
	\]
	Note that $A_{41} = A_{42} = 0$ by the Hessenberg condition, which corresponds to the case $(k,l) \notin H$.
	Then we can eliminate $g_{42}$ and $g_{41}$, using \eqref{eq:recursion}:
	\begin{align*}
		g_{42}
		&= \frac{1}{\lambda_4 - \lambda_2}g_{43}A_{32}
		= \frac{\lambda_3 - \lambda_2}{\lambda_4 - \lambda_2}g_{43}g_{32}, \\
		g_{41}
		&= \frac{1}{\lambda_4 - \lambda_1} ( g_{42}A_{21} + g_{43}A_{31} )\\
		&= \frac{\lambda_3 - \lambda_1}{\lambda_4 - \lambda_1}g_{43}g_{31} + \frac{(\lambda_2 - \lambda_1)(\lambda_3 - \lambda_4)}{(\lambda_4 - \lambda_1)(\lambda_4 - \lambda_2)}g_{43}g_{32}g_{21}.
	\end{align*}
\end{example}

Now we give a system of coordinates on the chart $X \cap U$.
\begin{lemma}
\label{lem:chart}
The functions $z_{kl}:=g_{kl}$, $(k,l)\in H$, identify $X\cap U$ with $\mathbb C^{|H|}$. More precisely,
\begin{equation}\label{eq:ring}
	\mathcal O(X\cap U) \cong \C[z_{kl} \mid (k,l)\in H].
\end{equation}
The diagonal torus preserves this chart and acts on its points by $z_{kl}(t\cdot q)=(t_k/t_l)z_{kl}(q)$. For $(i,j)\in H$, put $I_{ij}=\{(k,l)\in H \mid j\le l<k\le i\}$. Then
\begin{equation}\label{eq:subspace}
	X_{ij} \cap U = \{z_{kl} = 0 \mid (k,l)\notin I_{ij}\} \cong \mathbb C^{|I_{ij}|},
\end{equation}
and restriction of regular functions is the ring map
\begin{equation}\label{eq:restriction}
	\C[z_{kl} \mid (k,l)\in H] \to \C[z_{kl} \mid (k,l)\in I_{ij}],\qquad
	z_{kl} \mapsto
	\begin{cases}
		z_{kl}&(k,l)\in I_{ij},\\
		0&(k,l)\notin I_{ij}.
	\end{cases}
\end{equation}
\end{lemma}
\begin{proof}
Choose the entries $g_{kl}$ with $(k,l)\in H$ freely.
Induct on the height $k-l$.
When $(k,l)\notin H$ and $k>l$, we have $k>h(l)$, so $A_{kl}=0$ by \eqref{eq:hess}.
Equation \eqref{eq:recursion} then uniquely determines $g_{kl}$ as a polynomial in the free coordinates,
since $\lambda_k-\lambda_l\ne0$, proving \eqref{eq:ring}.

Diagonal conjugation sends $g$ to $tgt^{-1}$, so the stated formula for $z_{kl}$ follows.
The character $t_k/t_l$ restricts trivially to $\mathbb T_{ij}$ exactly when either $k=l$ or both $k,l$ belong to $[j,i]$:
outside that interval, the remaining diagonal characters are independent modulo the scalar diagonal.
Since our chart is an affine space with diagonal $\mathbb T_{ij}$-action, its fixed locus is the coordinate subspace in \eqref{eq:subspace}.
It is connected and contains $e$; hence it lies in the component $X_{ij}$, and the intersection of that component with the chart is exactly this fixed coordinate subspace. Formula \eqref{eq:restriction} follows.
\end{proof}

\begin{lemma}\label{lem:additive}
If a subgroup $\U$ of $C_G(\mathbb T)$ is isomorphic to the additive group $\mathbb{G}_a$, then $\U$ acts trivially on every $X_{ij}$.
\end{lemma}
\begin{proof}
The connected groups $\U$ and $\mathbb T$ preserve each connected component $X_{ij}$ of $X^{\mathbb T_{ij}}$.
By Lemma~\ref{lem:smaller}, $\Aut^0(X_{ij})\cong\operatorname{PGL}(V_{ij})$ (\cite[Theorems 2 and 3 on pp.\ 75--76]{akhi95}).
The quotient $\mathbb T/\mathbb T_{ij}$ is a maximal torus therein, whose centralizer is the torus itself.
As $\U$ centralizes $\mathbb T$, its image in $\Aut^0(X_{ij})$ lies in this torus.
Every algebraic-group homomorphism from $\mathbb{G}_a$ to a torus is trivial.
\end{proof}

\section{Nonexistence of centralizing additive subgroups}
\begin{lemma}
\label{lem:detection}
Let $\xi \in H^0(X, TX)^{\mathbb{T}}$ be a $\mathbb{T}$-invariant global vector field on $X$.
Its restriction $\xi|_{X_{ij}}$, which is a priori a section of $TX|_{X_{ij}}$, is in fact a vector field on each $X_{ij}$. 
Furthermore, if $\xi \neq 0$, then $\xi|_{X_{ij}} \neq 0$ for some $(i,j) \in H$.
\end{lemma}
\begin{proof}
Since $\xi$ is $\mathbb{T}_{ij}$-invariant, its value at any $q\in X_{ij}$ belongs to $(T_qX)^{\mathbb{T}_{ij}}=T_qX_{ij}$. Thus $\xi$ is tangent to every $X_{ij}$.

Work on the affine chart $X\cap U$ of Lemma~\ref{lem:chart}. Its tangent sheaf is free with basis $\partial/\partial z_{ij}$, so
\[
	\xi = \sum_{(i,j)\in H} f_{ij}(z) \frac{\partial}{\partial z_{ij}},\qquad
	f_{ij}\in \C[z_{kl} \mid (k,l)\in H].
\]
Since $X$ is irreducible and $\xi$ is a nonzero global section of $TX$, its restriction to the nonempty open chart $X\cap U$ is nonzero.
Hence at least one coefficient $f_{ij}$ is nonzero; fix such a pair $(i,j)$.
We consider the weights of the monomials occurring in $f_{ij}$.

Write $\mathbb T$-characters multiplicatively as $\chi_{kl}(t)=t_k/t_l$, while using $t_k-t_l$ for their additive notation. On points of the opposite big cell, diagonal conjugation sends the root coordinate to
\[
	z_{kl}(t\cdot q) = \chi_{kl}(t)z_{kl}(q).
\]
We use the induced action on functions $(t\cdot f)(q)=f(t^{-1}\cdot q)$.
Thus the coordinate function $z_{kl}$ is a $\mathbb T$-eigenvector of additive weight $t_l-t_k$, and the coordinate derivation $\partial/\partial z_{ij}$ has the opposite weight $t_i-t_j$.
Since the derivations form a free eigenbasis and $\xi$ is invariant, comparison of the coefficients of each derivation gives
\[
	t\cdot f_{ij} = \chi_{ji}(t) f_{ij},\text{ that is, }\wt(f_{ij}) = t_j-t_i.
\]
Explicitly, the polynomial $f_{ij}$ is of the following form:
\[
	f_{ij} = \sum_{m \in \mathbb Z_{\ge0}^{H}} c_m z^m,\qquad
	z^m = \prod_{(k,l)\in H} z_{kl}^{m_{kl}}.
\]
Each monomial is itself a $\mathbb T$-eigenvector, with weight
\[
	\wt(z^m)=\sum_{(k,l)\in H}m_{kl}(t_l-t_k).
\]
Since $f_{ij}$ is a $\mathbb T$-eigenvector of weight $t_j-t_i$,
every monomial $z^m$ occurring in $f_{ij}$ has weight $t_j-t_i$.
Hence, whenever $c_m\ne0$,
\begin{equation}\label{eq:weights}
\sum_{(k,l)\in H}m_{kl}(t_l-t_k)=t_j-t_i.
\end{equation}

To analyze the monomials in $f_{ij}$, we express the weights in terms of the simple roots $\alpha_k = t_k-t_{k+1}$.
Since $l<k$,
\[
	t_l - t_k = \alpha_l + \cdots + \alpha_{k-1},\qquad
	t_j - t_i = \alpha_j + \cdots + \alpha_{i-1}.
\]
The coefficients $m_{kl}$ are nonnegative.
Consequently, no term with $m_{kl}>0$ can contribute a simple root outside the interval $\{\alpha_j,\ldots,\alpha_{i-1}\}$, because such a contribution could not be cancelled by any other term.
Therefore
\[
	m_{kl}>0\quad\Rightarrow\quad j\le l<k\le i.
\]
Equivalently, $f_{ij}\in\C[z_{kl} \mid j\le l<k\le i]$.
By the explicit restriction map \eqref{eq:restriction}, $f_{ij}\vert{}_{X_{ij}}$ is nonzero.
Since this is the coefficient of $\partial/\partial z_{ij}$ in $\xi\vert{}_{X_{ij}}$, we conclude that $\xi\vert{}_{X_{ij}}\ne0$.
\end{proof}

\begin{proposition}\label{prop:no-ga}
There is no subgroup $\U$ isomorphic to $\mathbb{G}_a$ in $C_G(\mathbb T)\cap G_u$.
\end{proposition}
\begin{proof}
Suppose such a subgroup exists. Choose $0\ne v\in\Lie(\U)$.
Since the action of $G$ on $X$ is effective in characteristic zero, the induced vector field $\xi_v$ is nonzero.
Because $\U$ centralizes $\mathbb T$, the vector field is $\mathbb T$-invariant.
Lemma~\ref{lem:detection} gives a pair $(i,j)\in H$ for which $\xi_v|_{X_{ij}}\ne0$.
But Lemma~\ref{lem:additive} says that $\U$ acts trivially on $X_{ij}$, so the restricted vector field is zero, a contradiction.
\end{proof}

\section{The centralizer of the diagonal torus}
We recall that, for $1 \le k \le n-1$, $\alpha_k=t_k-t_{k+1}$ are the simple roots of type $\mathrm{A}_{n-1}$ and $\varpi_k$ are the corresponding fundamental weights.
The root order is defined by $\mu\le\nu$ when $\nu-\mu$ is a nonnegative integral combination of simple roots.

\begin{proposition}\label{prop:centralizer}
If $h\ne(n,\ldots,n)$, then $C_G(\mathbb T)=G$.
\end{proposition}
\begin{proof}
Suppose $C_G(\mathbb T)\ne G$.
Since $G$ is connected and algebraic, the conjugation representation of $\mathbb T$ on $\Lie(G)$ has a nonzero weight $\chi$ with nonzero weight space.
Choose a corresponding nonzero global vector field $\xi$ of weight $\chi$.
Since the Weyl group $S_n$ acts transitively on the Weyl chambers,
we can pick $w\in S_n$ such that $\mu = -w^{-1}\chi$ is dominant (e.g.\ see \cite[\S 13.2]{Hum72}).

At the fixed point corresponding to $w$, use the shifted cell $wU$, represented by matrices $wg$ with $g$ lower unitriangular.
Applying Lemma~\ref{lem:chart} to the diagonal matrix $w^{-1}Sw$ gives an affine chart $X \cap wU \cong \mathbb C^{|H|}$ with coordinates $z_{kl}^{(w)}$ whose weights as functions are $t_{w(l)} - t_{w(k)}$.
Express
\[
	\xi = \sum_{(i,j)\in H}f_{ij}^{(w)}(z)\frac{\partial}{\partial z_{ij}^{(w)}}.
\]
The restriction of a nonzero global vector field to this nonempty open chart is nonzero, so some polynomial coefficient contains a nonzero monomial $\prod (z_{kl}^{(w)})^{m_{kl}}$.
Comparing its weight with $\chi$ gives
\[
	\chi = \sum_{(k,l)\in H}m_{kl}(t_{w(l)}-t_{w(k)})+(t_{w(i)}-t_{w(j)}).
\]
Applying $-w^{-1}$ yields
\[
	\mu = (t_j-t_i) - \sum_{(k,l)\in H} m_{kl}(t_l-t_k) \le t_j-t_i.
\]
Since $\mu$ is a nonzero dominant weight, its expression in terms of the fundamental weights $\varpi_1, \dots, \varpi_{n-1}$ has nonnegative coefficients.
Fundamental weights
\[
	\varpi_k
	= \sum_{i=1}^k t_i - \frac kn\sum_{i=1}^n t_i 
	=\sum_{i=1}^{k}\frac{i(n-k)}n\alpha_i +\sum_{i=k+1}^{n-1}\frac{k(n-i)}n\alpha_i
\]
are expressed as linear combinations of simple roots with strictly positive coefficients.
Therefore so is $\mu$, and hence so is $t_j - t_i$ by the inequality $\mu\le t_j-t_i$.
Hence $t_j - t_i$ must be the highest root $t_1-t_n$, that is, $(i,j)=(n,1)$.
As $(n,1)\in H$, we have $h(1)=n$, and monotonicity gives $h=(n,\ldots,n)$, a contradiction.
\end{proof}

\section{Main theorem}
\begin{theorem}\label{thm:main}
Assume $h(i) \ge i+1$ for $i<n$ and $h \ne (n,\ldots,n)$.
Let $S$ be an $n \times n$ diagonal matrix with distinct diagonal entries and $\mathbb T$ the effective diagonal torus acting on $X = \Hess(S,h)$.
Then
\[
	\Aut^0(X) = \mathbb T \cong (\C^*)^{n-1}.
\]
\end{theorem}
\begin{proof}
We use the previously established decomposition $G \cong \mathbb T \ltimes G_u$ from \cite{JangEtAl}.
Proposition~\ref{prop:centralizer} gives $C_G(\mathbb T) = G$, hence $G_u=C_G(\mathbb T)\cap G_u$. 
Since $G_u$ is connected and unipotent, it contains a subgroup isomorphic to $\mathbb{G}_a$ whenever it is nontrivial.
This contradicts Proposition~\ref{prop:no-ga}.
Therefore $G_u = 1$ and $G = \mathbb T$.
\end{proof}
\begin{remark}
Theorem $\ref{thm:main}$ determines the identity component but does not, by itself, determine
the entire automorphism group.  Since $\Aut^0(X)=\mathbb T$, its component group is
\[
	\pi_0(\Aut(X)) = \Aut(X)/{\Aut^0(X)} = \Aut(X)/\mathbb T.
\]
It follows from \cite{JangEtAl} that this group is finite.
In fact, \cite{JangEtAl} shows that the component group is isomorphic to a subgroup of the automorphism group of the associated GKM graph.
This GKM graph automorphism group is $S_n$, or $S_n \times \mathbb Z/2$ when the index set $H$ admits the symmetry $(i,j) \leftrightarrow (n+1-j,n+1-i)$,
that is, the symmetry of the Dynkin diagram of type $\mathrm{A}_{n-1}$.
In other words, the component group is constrained by the symmetries of the GKM graph.

Suppose $n \ge 3$ and
\[
	h = (2,3,\ldots,n,n).
\]
Then $X$ is the permutohedral toric variety, whose fan is the braid fan of type $\mathrm{A}_{n-1}$.
Demazure's general description of the automorphism group of a complete nonsingular toric variety (see \cite{Demazure,Oda}) applies directly in this case.
The braid fan has no Demazure roots, so the identity component is the acting torus $\mathbb T$.
Moreover, the component group is the automorphism group of the braid fan.
This fan is preserved by permutations of the $n$ coordinates and by the central symmetry $v\mapsto-v$;
these generate its full automorphism group.
Consequently,
\[
	\pi_0\bigl(\Aut(\Hess(S,h))\bigr) \cong S_n \times\mathbb Z/2.
\]
Thus this special case gives the largest component group allowed by the general
description above.

When $h = (n-1,n,\ldots,n)$, or equivalently, when $X$ has codimension one in $\Fl(\C^n)$,
a more precise description is given in \cite{BEHLLMS} for $n\ge 4$: the component group is expressed
in terms of the stabilizer of the eigenvalue configuration, together with the involution induced by the symmetry of the Dynkin diagram.
In particular, the component group may vary with the unordered eigenvalue configuration of $S$,
considered up to affine transformations.
\end{remark}

\section*{Acknowledgements}
The proof presented in this paper was suggested to the authors through discussions with ChatGPT, developed by OpenAI.
The authors subsequently worked out and verified the mathematical details and take full responsibility for the contents of the paper.

The first author is supported in part by JSPS KAKENHI Grant-in-Aid for Scientific Research (C): 25K07007.
The second author is supported in part by JSPS KAKENHI Grant-in-Aid for Scientific Research (C): 25K00205.
The third author is supported in part by NSFC: 11901218.

\end{document}